\documentclass[draft,12pt,leqno]{amsart}

\usepackage{amsmath,amssymb,amscd,amsthm}
\usepackage{latexsym}
\usepackage[dvips]{graphics,epsfig}

\newtheorem{theorem}{Theorem}
\newtheorem{lemma}{Lemma}
\newtheorem{proposition}{Proposition}

\newtheorem{corollary}{Corollary}

\newtheorem{claim}{Claim}

\newcommand{\q}{\quad}
\newcommand{\qq}{\quad\quad}

\newcommand{\norm}[2]{{\left\| #1 \right\|}_{#2}}
\newcommand{\f}[2]{\frac{#1}{#2}}

\newcommand{\al}{\alpha}

\newcommand{\de}{\delta}

\newcommand{\De}{\Delta}
\newcommand{\ve}{\varepsilon}

\newcommand{\si}{\sigma}

\newcommand{\rd}{{\mathbf R}^d}

\newcommand{\rone}{\mathbf R^1}
\newcommand{\rtwo}{\mathbf R^2}

\newcommand{\cz}{\mathcal Z}

\newcommand{\cc}{\mathcal C}

\newcommand{\intl}{\int\limits}

\newcommand{\suml}{\sum\limits}

\newcommand{\supl}{\sup\limits}

\newcommand{\p}{\partial}

\newcommand{\beq}{\begin{equation}}
\newcommand{\eeq}{\end{equation}}
\newcommand{\beqna}{\begin{eqnarray*}}
\newcommand{\eeqna}{\end{eqnarray*}}
\newcommand{\beqn}{\begin{equation*}}
\newcommand{\eeqn}{\end{equation*}}
\newcommand{\bp}{\begin{proof}}
\newcommand{\ep}{\end{proof}}
\newcommand{\bprop}{\begin{proposition}}
\newcommand{\eprop}{\end{proposition}}
\newcommand{\bt}{\begin{theorem}}
\newcommand{\et}{\end{theorem}}
\newcommand{\bex}{\begin{Example}}
\newcommand{\eex}{\end{Example}}
\newcommand{\bc}{\begin{corollary}}
\newcommand{\ec}{\end{corollary}}
\newcommand{\bcl}{\begin{claim}}
\newcommand{\ecl}{\end{claim}}
\newcommand{\bl}{\begin{lemma}}
\newcommand{\el}{\end{lemma}}

\begin{document}

\title
[Gevrey estimates for the Kuramoto-Sivashinsky equation]
{The Kuramoto-Sivashinsky equation in $R^1$ and  $R^2$: 
effective estimates of the 
high-frequency tails and higher Sobolev norms}

\author{Milena Stanislavova} 
\author{Atanas Stefanov}

\address{Milena Stanislavova\\
Department of Mathematics \\
University of Kansas\\
1460 Jayhawk Blvd\\ Lawrence, KS 66045--7523}
\email{stanis@math.ku.edu}
\address{Atanas Stefanov\\
Department of Mathematics \\
University of Kansas\\
1460 Jayhawk Blvd\\ Lawrence, KS 66045--7523}

\email{stefanov@math.ku.edu}

\thanks{
Stanislavova's research supported in part by  
NSF-DMS 0508184. 
Stefanov's research supported in part by  
NSF-DMS 0701802.}
\date{\today}



\begin{abstract}
We consider the Kuramoto-Sivashinsky (KS) equation in finite domains 
of the form
$[-L,L]^d$. Our main result provides
refined Gevrey estimates for the solutions of the one dimensional differentiated KS, which
in turn imply effective new estimates for higher Sobolev norms of the 
solutions in terms of powers of $L$. We illustrate our method on a 
simpler model,
namely the regularized Burger's equation. We also show local 
well-posedness for 
the two dimensional KS  equation and provide an explicit criteria
for (eventual) blow-up in terms of its $L^2$ norm. 
The common underlying idea in both results is that {\it a priori} 
control of the $L^2$ norm  is enough in order to conclude higher order
regularity and allows one  to get good estimates on the 
high-frequency tails of the solutions. 
\end{abstract}

\maketitle

\section{Introduction}
The Kuramoto-Sivashinsky equation 
\begin{equation}
\label{eq:1}
\left\{\begin{array}{l} 
\phi_t +\De^2 \phi+\De \phi + \f{1}{2}|\nabla\phi|^2 =0 \qq x\in [-L,L]^d \\
\phi(t,x+ 2 Le_j)=\phi(t,x), j=1, \ldots, d. \\
\phi(0,x):=\phi_0(x)
\end{array} \right.
\end{equation}
where $d\geq 1$ and $L>0$ models pattern formation in different physical contexts. It arises as a model of nonlinear evolution of linearly unstable interfaces in a variety of applications such as flame propagation (advocated by Sivashinsky \cite{Siv}) and reaction-difusion systems (derived by Kuramoto in \cite{Ku}). It has been studied extensively by many authors. It is interesting mathematically because the linearization about the zero state has a large number of exponentially growing modes. In \cite{wstrauss}, the instability of the travelling waves is a hint of the complexity of the dynamics of KS equation in the unbounded case. The main results in the periodic case are on the global existence of the solutions, their stability and long-time behavior.

In one space dimension, it is convenient to consider  
the differentiated Kuramoto - Sivashinsky
equation. That is, set $u=\phi_x$ and differentiate \eqref{eq:1} with respect to
$x$ to get a closed form equation for $u$
\begin{equation}
\label{eq:2}
\left\{\begin{array}{l} 
u_t + u_{xxxx}+u_{xx} + u u_x  =0   \\
u(t,x+ 2L)=u(t,x) \\
u(0,x):=u_0(x)
\end{array} \right.
\end{equation}
In this periodic case, the global well-posedness of \eqref{eq:2}, 
the existence of global attractor and it's dimension were studied in 
\cite{BG, GO, NST, Collet} and many others. Of interest here is the existence of attracting ball and the dependence
of  $\lim_{t\to \infty} \|u(t, \cdot)\|_{L^2}$ on the size of the domain $L$. 
The best possible current result 
$\lim_{t\to \infty} \|u(t, \cdot)\|_{L^2}=o(L^{3/2})$ is achieved by Giacomelli and
Otto in \cite{GO}, see also \cite{BG} for a somewhat more direct  proof of 
the slightly weaker result 
$\lim_{t\to \infty} \|u(t, \cdot)\|_{L^2}\lesssim L^{3/2}$. 
We would like to point out
that this last bound applies as well to the solutions of  the
so-called destabilzied KS equation 
\begin{equation}
\label{eq:45}
u_t + u_{xxxx}+u_{xx} - \eta u + u u_x  =0, \q\eta>0
\end{equation}
and moreover, such result is {\it optimal} in this context. 
Using techniques similar to \cite{BG}, the authors of 
 \cite{BFG} consider a  nonlocal Kuramoto-Sivashinsky equation  and prove 
estimates for $\lim_{t\to \infty} \|u(t, \cdot)\|_{L^2}$. 
In this case one gets different estimate in the case of odd initial data 
from the case  of arbitrary initial data. 

Before we embark on our discussion on the optimality of these results, it is worth 
noting the following two conjectures. Namely, based on numerical simulations 
about the dimension of the attractor,  it is conjectured that 
$\|u(t, \cdot)\|_{L^2}$ behaves according to 
\begin{equation}
\label{eq:l2}
\limsup_{t\to \infty} \|u(t, \cdot)\|_{L^2}\leq C L,
\end{equation}
whereas for $\|u(t, \cdot)\|_{L^\infty}$ 
\begin{equation}
\label{eq:011}
\limsup_{t\to \infty} \|u(t, \cdot)\|_{L^\infty}\leq C.
\end{equation}
If true, these would be the best possible estimates, since these are satisfied by the stationary
solutions of the problem, see \cite{Michelson}. 
For a nice discussion about these conjectures 
the reader is referred to the introduction in \cite{CT}.

In two space dimensions, even the question of global regularity of the Cauchy problem 
$$
\phi_t +\De^2 \phi+\De \phi + \f{1}{2}|\nabla\phi|^2 =0, \qq  \phi(0,x):=\phi_0(x)
$$
in $R^N, N\geq 2$ or in the periodic boundary conditions case is still open. The results in 
\cite{SellTaboada} and \cite{Molinet} show local existence and local dissipativity 
with some restrictions on the domain and the initial data. 
In this direction the best result so far is in \cite{BG}, showing that in $L^2((0,L_x)\times(0,L_y))$ with $L_y\leq CL_x^{13/7}$ one has
$\displaystyle \limsup_{t\rightarrow \infty}
\norm{\phi} \leq CL_x^{3/2}L_y^{1/2}.$  In the present work, we will show that 
the solution is defined and classical up to time $T^*\leq \infty$, provided 
$\limsup_{t\to T^*}\|\phi(t, \cdot)\|_{L^2}<\infty$. In fact, we will be able to present 
an explicit Gronwall's type argument, which allows one to control higher Sobolev norms so long
as $\|\phi(t, \cdot)\|_{L^2}$ is under control.

The question of Gevrey class regularity for the Kuramoto-Sivashinsky equation is of interest because it can be used to improve the error estimates in the computation of the approximate inertial manifolds (see \cite{JKT} and also \cite{FT},\cite{Titi} 
for a similar results on the Navier-Stokes equations). 
In \cite{Liu} the author studies the Gevrey class regularity 
for the odd solutions of the one dimensional Kuramoto-Sivashinsky 
equation with periodic boundary conditions and odd initial data. Theorem 1 
in his paper should be compared with the estimates in Corollary \ref{cor:1} 
and Corollary \ref{cor:2} of the current paper, see the remarks after Corollary \ref{cor:lo}.

Our main results are Gevrey regularity theorems for the 
solutions of \eqref{eq:2}, but we will not emphasize our 
presentation on that fact. Instead, we will concentrate 
on the specific estimates that one can get for the high-frequency
tails of the solutions of \eqref{eq:2}. 

In order to illustrate our methods on a somewhat simpler model, 
we will first consider the regularized Burger's equation. In it, we can actually take the regularization operator in the form 
$A_s=(-\De)^{s/2}$.Thus, our model is 
\begin{equation}
\label{eq:30}
\left\{\begin{array}{l}  
u_t - A_s u +div(u^2)=0\qq x\in [-L,L]^d \\
u(t,x+ 2L)=u(t,x) \\
u(0,x):=u_0(x),
\end{array} \right.
\end{equation}
where the  formal definition of $A_s$ is given in Section \ref{sec:four}.  

In the next two theorems, we give estimates of the high-frequency tails  of the
solutions of \eqref{eq:30} and \eqref{eq:2} respectively. For this, 
we shall need 
the Littlewood-Paley projections, which are defined in Section 
\ref{sec:four} below.
We have 
\begin{theorem}
\label{theo:2}
Let $d\geq 1$, $1<s\leq 2$ or $s>1+d/2$. Then, the regularized Burger's equation 
\eqref{eq:30} is a globally well-posed problem, whenever
the data belongs to $L^2$. \\ In
addition, in the case $1<s\leq 2$, assume  $u_0\in L^2\cap L^\infty$. 
Then, for every $1>>\de>0$, there exists $C_{\de,s}$, so that for any  $j\geq 0$,
\begin{equation}
\label{eq:31}
\|P_{\geq 2^j L} u(t, \cdot)\|_{L^2}^2\leq  
(C_{\de,s}  \max(1, \|u_0\|_{L^2\cap L^\infty}^2))^{j+1} 2^{-\min(t,1)
(1-\de)(s-1)j^2}. 
\end{equation}
For $s>1+d/2$,  one has a constant $C_{s}$
\begin{equation}
\label{eq:310}
 \|u_{>2^j L} (t, \cdot)\|_{L^2}^2\leq 
 (C_{\de,s}  \max(1, \|u_0\|_{L^2}^2))^{j+1} 2^{-\min(t,1)(s-1-d/2)j^2}.  
\end{equation}
As an easy  corollary, one can 
estimate  $ \sup_{\de<t<\infty}\|u(t, \cdot)\|_{H^m} $ in terms of quantities, 
which are independent of the size of the domain $L$. 
\end{theorem}
{\bf Remarks:} 
\begin{itemize}
\item In both cases, our results show that 
 the solution belongs to the Gevrey class $\mathcal{G}^2$. In particular the
 function $x\to u(t, x)$ is real-analytic for every fixed $t>0$.
\item The results in Theorem \ref{theo:2} can be extended accordingly 
to the case of $\rd$. 
\end{itemize}
Similar results hold for the one dimensional Kuramoto-Sivashinsky equation
 \eqref{eq:2}. The main difference with the regularized Burger's equation will be the 
 unavailability of control of $\|u(t, \cdot)\|_{L^2}$ over the course of the evolution. 
 In fact, as  discussed previously, the function $t\to \|u(t, \cdot)\|_{L^2}$ may 
 (and sometimes does) grow to at least $C \sqrt{L}$ for the static solutions of \eqref{eq:2}, 
 see \eqref{eq:l2}. 
 \newpage
\begin{theorem}
\label{theo:3}
Let  $u_0\in L^2(-L,L)$ and $L>>1$. Set $H=\sup_{0\leq s<\infty} 
\|u(s)\|_{L^2}$, where $u$ is a solution to \eqref{eq:2}. 
 Then, there exist  absolute constants $C_0, C_1$, so
that for every $j\geq 0$, 
\begin{equation}
\label{eq:42}
\|u_{>C_0 2^j H^{2/5}L}(t, \cdot)\|_{L^2}\leq 
C_1^j 2^{-\f{1}{2}\min(t,5/2) j^2}\sup_s\|u(s)\|_{L^2}.
\end{equation}
\end{theorem}
Regarding $H$ in the statement in Theorem \ref{theo:3}, 
one may actually infer from the results of \cite{BG}, \cite{GO} and the earlier
papers on the subject that 
\begin{equation}
\label{eq:010}
H=\supl_s \|u(s, \cdot)\|_{L^2}\leq \|u_0\|_{L^2}o(1/t) + C L^{3/2},
\end{equation}
of which then $\limsup_{t\to \infty}\|u(t, \cdot)\|_{L^2}=O(L^{3/2})$ 
is a corollary. Thus, when $L>>1$ (in particular when $\|u_0\|_{L^2}<<L$), 
we have
that $H\leq C L^{3/2}$. In particular, we have an estimate of 
$\|u_{\gtrsim L^{8/5}}(t, \cdot)\|_{L^2}$, but we prefer to formulate this 
as  estimates on the higher
Sobolev norms. 
\begin{corollary}
\label{cor:lo}
Let  $s\geq 0$, $L>>1, \de>0$ and $\|u_0\|_{L^2}<<L$.  
Then, there exists $C_{s, \de}$, so that 
\begin{equation}
\label{eq:43}
\supl_{\de\leq t<\infty} \|u(t, \cdot)\|_{H^s}\leq C_{s, \de} L^{3s/5} L^{3/2}.
\end{equation}
\label{cor:1}
\end{corollary}
{\bf Remarks on Theorem \ref{theo:3} and Corollary \ref{cor:lo}}
\begin{enumerate}
\item The estimate \eqref{eq:43} may be stated (with the same assumptions on
$\|u_0\|_{L^2}$) 
in the form 
\begin{equation}
\label{eq:47}
\supl_{\de\leq t<\infty} \|u(t, \cdot)\|_{H^s}\leq C_{s, \de}
H^{2s/5+1}.
\end{equation}
In other words, if one improves the bounds on $H$ in \eqref{eq:010}, 
then one immediately gets an improvement of the results \eqref{eq:43}
in the form \eqref{eq:47}. Said  differently, 
 with the best current technology, namely
$H\lesssim L^{3/2}$, \eqref{eq:43} is an instance of \eqref{eq:47}. 
\item 
The bounds \eqref{eq:43} and \eqref{eq:47} apply for solutions of the 
destabilized Kuramoto - Sivashinsky
equation \eqref{eq:45} as well. As we have discussed previously, $H\lesssim
L^{3/2}$ is optimal here in contrast with the standard KSE. 
\item The estimate \eqref{eq:43} should be compared with the bound on 
$\sup_t \|u(t,\cdot)\|_{H^s}$ by Liu,  \cite{Liu}, which is of the form 
$\sup_t \|u(t,\cdot)\|_{H^s}\lesssim L^{4s+5/2}$ and which follows from a similar Gevrey
regularity estimate.  One should have in mind that the best available bound at the
time\footnote{which Liu has used in his estimates} was 
$\sup_t \|u(t,\cdot)\|_{L^2}\lesssim L^{5/2}$. Even with the use of that 
 bound however, our method from Theorem \ref{theo:2} would have produced an estimate of the form  
 $\sup_t \|u(t,\cdot)\|_{H^s}\lesssim L^{s+5/2}$, which is again superior to the results of 
 \cite{Liu}.
\end{enumerate}
Next, we present another estimate, which gives bounds on $\sup_t \|u(t,\cdot)\|_{H^s}$ in 
terms of $\sup_t \|u(t, \cdot)\|_{L^\infty}$. This follows essentially the same scheme of proof
and yet, it gives at least as good bounds\footnote{and potentially much better bounds} as \eqref{eq:43}, see the discussion after 
Corollary \ref{cor:2}. The reason for the effectiveness of such an approach is that 
it almost avoids the use of Sobolev embedding, which is clearly ineffective in this
context.

 It is actually possible to give yet another different form of the estimates in 
Corollary \ref{cor:1} in terms of the 
quantities\footnote{As it was pointed out already,
there is the standing conjecture  \eqref{eq:011}, 
which  puts an  uniform bound
on $K_\infty$. } 
$K_p=\sup_{0<t<\infty} \|u(t, \cdot)\|_{L^p}$, where one should think of 
$p$ as being very large.    
\begin{corollary}
\label{cor:2} Let $s\geq 0$. Then, there exists a constant 
$C_{s,p,\de}$, so that 
\begin{equation}
\label{eq:012}
\supl_{\de\leq t<\infty} \|u(t, \cdot)\|_{H^s}\leq C_{s,p, \de}
K_p^{s/(3-1/p)} H.
\end{equation}
Roughly speaking, we get a factor of $K_\infty^{1/3}$ for every derivative of $u$. 
\end{corollary}
{\bf Remark:} We would like to point out that to the best of our knowledge, 
 the best   estimate  currently 
available for $K$ is obtainable through the Sobolev embedding theorem and the estimates for 
 $\sup_t \|u(t)\|_{H^{1/2+}}$ from \eqref{eq:47}. 
 This is  certainly a  very crude estimate, but let us  use it anyways.  
By the bound $H\lesssim L^{3/2}$ and assuming
 $\|u_0\|_{L^2}<<L$, $L>>1$,  we have that  for every $2<p<\infty$
 by\footnote{Here we are ignoring the minor issue for the bounds in the interval
 $0<t<\de$, but recall that our discussion is about global behavior.}
 \eqref{eq:47}
 $$
 K_p\leq C_p\sup\|u(t, \cdot)\|_{H^{1/2-1/p}}\leq C_p H^{1/5+1-2/(5p)}\leq C_p 
 L^{9/5-2/(5p)}.
 $$
Clearly, with this bound for $K_p$, \eqref{eq:012} is only slightly worse than 
\eqref{eq:43}. However, if the conjecture \eqref{eq:011} holds true or even an estimate of the
form $K_\infty\lesssim L^{9/5-}$ is established, then
\eqref{eq:012} gives better result. Indeed, if \eqref{eq:011} holds, then  
\begin{equation}
\label{eq:013}
\supl_{\de<t<\infty} \|u(t, \cdot)\|_{H^s}\leq C_{s,\ve, \de}
L^{\ve s} H
\end{equation}
 for every $\ve>0$. This would  one more time confirm the empirical observations,  that the 
 whole action in  the evolution of the KS comes in the low frequencies.

The following result concerns solutions for the KS 
equation \eqref{eq:1} in two spatial dimensions. More specifically, 
it characterizes the (eventual) blow-up time. 
\begin{theorem}
\label{theo:1}
Let $d=2$. Then, the KS equation \eqref{eq:1} is locally well-posed in the
following sense - for every initial data $\phi_0\in L^2([-L,L]^2)$, there exists
a time $T^*=T^*(\|\phi_0\|_{L^2})$, so that \eqref{eq:1} has an
unique classical solution $\phi$, 
$\phi(t, \cdot)\in C^\infty(\rtwo)\cap L^2(\rtwo)$ up to time $T^*$. 
In addition, the solution is either global (i.e. $T^*=\infty$) or else, it must
be that  
\begin{eqnarray*}
& & 
\lim_{t\to T^*-} \|\phi(t,\cdot)\|_{L^2}= \infty \\
& & 
 \lim_{t\to T^*-} \int_0^t 
\|\nabla \phi(t,\cdot)\|_{L^2}^2 dt=\infty.
\end{eqnarray*}
In other words, the solution is well-defined and classical up to time 
$T$ as long as either \\
$\lim_{t\to T-} \|\phi(t,\cdot)\|_{L^2}<\infty$ or $\lim_{t\to T-} \int_0^t 
\|\nabla \phi(t,\cdot)\|_{L^2}^2 dt<\infty$.
\end{theorem}
We would like to point out that the same theorem
applies in the case of three spatial dimensions. Its proof however requires an
additional step and we do not pursue it for the sake of brevity.

\section{Preliminaries} 
\label{sec:four}
Since our attention will be focused on the case of  domains 
$[-L,L]^d$, we will briefly introduce some relevant concepts from Fourier
series, which will be useful in the sequel. 
\subsection{Discrete Fourier transform and Plancherel's identity}
On the interval  $[-L,L]$, introduce the  Fouier transform 
$L^2([-L,L])\to l^2(\cz^d)$, 
by setting $f\to \{a_k\}_{k\in\cz^d}$, where 
$$
a_k= (2L)^{-d/2} \intl_{[-L,L]^d}  f(x) e^{-2\pi i k\cdot x/L} dx.
$$
The inverse Fourier transform  is the 
familiar Fourier expansion 
\begin{equation}
\label{eq:015}
f(x)=\f{1}{(2L)^{d/2}} \suml_{k \in \cz^d}  a_k e^{2\pi i k\cdot  x/L} .
\end{equation}
and the Plancherel's identity is 
$\norm{f}{L^2([-L,L]^d)}=\norm{\{a_k\}}{l^2(\cz^d)}$. 
Note that here and for the rest of the 
paper $L^2([-L,L]^d)$ is the space of 
square integrable functions with period $2L$ in all variables. 
In our case, we will be dealing
with real-valued functions only. 
\subsection{Littlewood-Paley projections and Bernstein inequality}
The Littlewood-Paley operators acting on $L^2([-L,L])$ are   defined for a 
function $f$ in the form of \eqref{eq:015} via 
$$
P_{\leq N} f (x) = \f{1}{(2L)^{d/2}} 
\suml_{k: |k|\leq N} a_k e^{2\pi i k\cdot  x/L}.
$$
That is $P_{\leq N}$ truncates the terms in the 
Fourier series expansion with frequencies 
$k: |k|> N$.  Clearly $P_{\leq N}$ is a projection 
operator. More generally, we may define for all $0\leq N<M\leq \infty$
$$
P_{N\leq \cdot\leq M}f(x) = \f{1}{(2L)^{d/2}} 
\suml_{k: N\leq |k|\leq M} a_k e^{2\pi i k\cdot  x/L}.
$$
Clearly, we may take $M, N$ to be nonintegers as well. A basic result in harmonic analysis 
on the torus is that Fourier series $P_{<N} f$ converge to $f$ in $L^p, 1<p<\infty$ sense. 
This is in fact equivalent to the uniform boundedness of the operators $P_{<N}$  in
$L^p([-L,L]^d)$, which we now record
$$
\|P_{<N} f\|_{L^p([-L,L]^d)}\leq C_{d,p} \|f\|_{L^p}.
$$
Note that this estimate fails as $p=\infty$ and thus $C_{s,p}\to \infty$ 
as $p\to \infty$. 

We will also need a  Sobolev embedding type result  for the spaces $L^q([-L,L]^d)$.  
We state it in the form of the {\it Bernstein inequality}.  
\begin{lemma}
\label{le:bern}
Let $N$ be an integer and $f:[-L,L]^d \cc$. Then, for 
every $1\leq p\leq 2 \leq q\leq \infty$, 
$$
\norm{P_{<N} f}{L^q}\leq C_{d,p,q} (N/L)^{d(1/p-1/q)} \norm{f}{L^p}.
$$
\end{lemma}
\begin{proof} 
The proof of this lemma is classical and can be found\footnote{in the case $L=1$, but
the general case follows easily by rescaling}, 
as Lemma 3  in \cite{StSt}.
\end{proof}
Next, we introduce the Sobolev spaces 
\begin{eqnarray*}
& & \dot{H}^s((-L,L)^d) = \{f:(-L,L)^d \to \cc\qq 
(\suml_{k\in \cz^d}|a_k|^2 \left(\f{|k|}{L}\right)^{2s} )^{1/2}<
\infty \}, \\
& & 
H^s  = L^2 \cap \dot{H}^s.
\end{eqnarray*}
One may also find convenient to work with the equivalent norm 
\begin{equation}
\label{eq:020}
\|f\|_{\dot{H}^s}\sim \left(\suml_{j\in \cz} 2^{2sj}  (\suml_{|k|\sim 2^j L} 
|a_k|^2) \right)^{1/2}\sim \left(\suml_{j\in \cz} 2^{2sj}  \|P_{\sim 2^j L} f\|_{L^2}^2 
 \right)^{1/2},
\end{equation}
which we will use regularly in the sequel. 
Another useful object to define is the (fractional) 
differentiation operator $A_s=(-\De)^{s/2}$, defined\footnote{The definition here
makes sense only for sequences $\{a_k\}$ with enough decay, 
say in  $l^2_{\si}, \si>s+d/2$. One may of course take $A_s f$ to 
represent a distribution for
less decaying $\{a_k\}$.}
  via 
$$
A_s [\suml_k a_k e^{2\pi i k \cdot x /L}] = \suml_k a_k 
\left(\f{2\pi |k|}{L}\right)^se^{2\pi i k \cdot x/L}.
$$
Sometimes in the sequel, we will just use the notation $|\nabla|^s$ instead of $A_s$. 
An useful  corollary of the representation \eqref{eq:020} is 
$$
\|A_s P_{2^j L} f\|_{L^2}\sim 2^{js} \|P_{2^j L} f\|_{L^2},
$$
and its obvious generalization $\|A_s P_{>2^j L} f\|_{L^2}\gtrsim 2^{js} 
\|P_{2^j L} f\|_{L^2}$ for $s\geq 0$. \\
The following simple orthogonality lemma is used frequently in the energy estimates
presented below.
\begin{lemma}
\label{le:90}
Let $A, B, C$ are three subsets of $\cz^d$, so that $0\notin A+B+C$. Then, for
any three functions $f,g,h\in L^2([-L,L]^d)$, 
\begin{equation}
\label{eq:016}
\intl_{[-L,L]^d} (P_A f)(P_B g)(P_C h) dx=0.
\end{equation}
As an useful corollary, for every $N$, 
\begin{equation}
\label{eq:017}
\intl_{[-L,L]^d} f_{>N} g_{<N/2} h dx= \intl_{[-L,L]^d} f_{>N} g_{<N/2} 
h_{>N/2} dx
\end{equation}
\end{lemma}
\begin{proof}
The proof of \eqref{eq:016} follows by expanding in Fourier series
$$
f g h(x)=(2L)^{-d/2}\suml_{k,m,n} f_k g_m h_n e^{2\pi i (k+m+n)\cdot x/L},
$$
and then realizing that since $(k+m+n)\neq 0$, all the terms will upon integration 
in $x$ result in zero. 
The proof of \eqref{eq:017} follows by observing that the difference between the
two sides is 
$$
\intl_{[-L,L]^d} f_{>N} g_{<N/2} 
h_{\leq N/2} dx=0,
$$
by \eqref{eq:016}, since $0\notin \{n:|n|>N\}+\{m:|m|<N/2\}+\{k:|k|\leq N/2\}$. 
\end{proof}

\section{Estimates of the high-frequency tails for regularized Burger's equations}
\label{sec:10}
In this section, we show that Theorem \ref{theo:2} holds. As we have pointed out
already, the essence of this result is a Gevrey  regularity of 
 the solution.
The classical theory guarantees global 
existence of classical solutions, so we proceed with the estimates. 
 
For $M>>L$, so that $M/L\in 2^{\cz}$,  
take the projection $P_{>M}$ on both sides of
\eqref{eq:30}.
We then take a scalar product of the result with $u$. We have 
\begin{equation}
\label{eq:35}
\f{1}{2} \p_t  \|u_{>M} (t, \cdot)\|_{L^2}^2+ \| P_{>M} 
A_s^{1/2} u(t, \cdot)\|_{L^2}^2\leq
|\int u_{>M} div(u^2) dx|
\end{equation}
Clearly, 
$$
\| P_{>M} A_s^{1/2} u(t, \cdot)\|_{L^2}^2\geq (M/L)^{s} \|u_{>M}(t, \cdot)\|_{L^2}^2, 
$$
while since $\int u_{>M} div[(u_>M)^2] dx=\f{1}{3} \int div[(u_>M)^3] dx=0$, one has 
\begin{eqnarray*}
& & 
\int u_{>M} div(u^2) dx= 2\int u_{>M} div[u_{>M} u_{\leq M}]dx+
\int u_{>M} div[ u_{\leq M}^2]dx= \\
& & = -2\int div(u_{>M}) u_{>M} u_{\leq M}+ 2 \int u_{>M} u_{\leq M} 
div[u_{\leq M}]dx = \\
& &=   \int u_{>M}^2 div(u_{\leq M})dx+ 2 \int u_{>M} u_{\leq M} 
div[ u_{\leq M}]dx\leq \\ 
& & \leq  \|u_{>M}\|_{L^2}^2
\|\nabla u_{\leq M}\|_{L^\infty}+ 2 \int u_{>M} u_{\leq M} 
div[ u_{\leq M}]dx,
\end{eqnarray*}
Furthermore, by Lemma \ref{le:90}
$\int u_{>M} u_{<M/2} 
div[ u_{<M/2}]dx=0$ and hence
\begin{eqnarray*}
& & 
\int u_{>M} u_{\leq M} 
div[ u_{\leq M}]dx= \int u_{>M} (u_{\leq M/2} +u_{M/2<\cdot \leq M}) 
div[u_{\leq M/2} +u_{M/2<\cdot \leq M}]=\\
& & = \int u_{>M} u_{\leq M/2} div[u_{M/2<\cdot \leq M}] dx+ 
\int u_{>M} u_{M/2<\cdot \leq M} div[u_{\leq M}]dx.
\end{eqnarray*}
The last identity allows us to estimate by H\"older's as follows
\begin{eqnarray*}
& & |\int u_{>M} u_{\leq M} 
div[ u_{\leq M}]dx|\leq C \|u_{>M}\|_{L^2} 
\|\nabla u_{M/2<\cdot\leq M}\|_{L^2} \|u_{\leq M/2}\|_{L^\infty}   \\
& & +C \|u_{>M}\|_{L^2} 
\|u_{>M/2}\|_{L^2} \|\nabla u_{\leq M}\|_{L^\infty} \leq  \\
& & \leq C(M/L) 
\|u_{>M}\|_{L^2}\|u_{>M/2}\|_{L^2} (\| u_{\leq M/2}\|_{L^\infty}\| + 
\|u_{\leq M}\|_{L^\infty})
\end{eqnarray*}
Inserting all the relevant estimates  in  \eqref{eq:35}  yields 
\begin{eqnarray*}
& & \p_t  \|u_{>M} (t, \cdot)\|_{L^2}^2+ 2 (M/L)^{s} 
\|u_{>M}(t, \cdot)\|_{L^2}^2\leq  \\
& & \leq C (M/L) \|u_{>M}\|_{L^2} \|u_{>M/2}\|_{L^2}
(\| u_{\leq M/2}\|_{L^\infty}\| + 
\|u_{\leq M}\|_{L^\infty}) 
\end{eqnarray*}
At this stage, the argument splits into the two cases, $1<s\leq 2$ and $s>1+d/2$. 
\subsection{Estimates in the case $1<s\leq 2$}
In this case, we 
use the results of \cite{Cordoba} (see also \cite{Ning}), where the authors have
established the following pointwise inequality\footnote{More precisely,
C\'ordoba-C\'ordoba established \eqref{eq:311} for $p=2^l, l=1, 2, \ldots$, while Ju,
\cite{Ning} has extended it in the range $2\leq p<\infty$.}
\begin{equation}
\label{eq:311}
\int_{[-L,L]^{d}} |\psi|^{p-2} \psi A^s[\psi] dx \geq C_{L,p} \|A^{s/2}
\psi^{p/2}\|_{L^2}^2.
\end{equation}
for any $0\leq s\leq 2$,  and for any 
smooth function $\psi(x):[-L, L]^d\to \rone$. Due to this inequality, 
one observes that
taking a scalar product of \eqref{eq:30} with $|u|^{p-2}u$ yields 
$$
\p_t \f{1}{p} \|u\|_{L^p}^p\leq \p_t \f{1}{p} \|u\|_{L^p}^p+ C_{L,p} \|A^{s/2}
\psi^{p/2}\|_{L^2}^2\leq  \int u_t  u |u|^{p-2}   dx + \int [A_s u] u |u|^{p-2} dx=0,
$$
whence $\|u(t, \cdot)\|_{L^p}$ is a decreasing function for every $p\geq 2$.\\
By Lemma \ref{le:bern} and 
the monotonicity of $t\to \|u(t, \cdot)\|_{L^p}: 2\leq p<\infty$, 
we have\footnote{This additional step is required , since the Littlewood-Paley
operators $P_{<M}$ are not bounded on $L^\infty$, otherwise, we would have preferred
to take $p=\infty$ and not lose the factor $(M/L)^{d/p}$.}
\begin{equation}
\label{eq:38}
\| u_{\leq M/2}\|_{L^\infty}\| + 
\|u_{\leq M}\|_{L^\infty}\leq C_{d,p} (M/L)^{d/p}\|u_0\|_{L^p}.
\end{equation}
for any $p: 2<p<\infty$. Select $p: d/p=2\de(s-1)$.  Thus, after
Cauchy-Schwartz's inequality 
\begin{equation}
\label{eq:90}
\p_t  \|u_{>M} (t, \cdot)\|_{L^2}^2+ (\f{M}{L})^{s} 
\|u_{>M}(t, \cdot)\|_{L^2}^2\leq C_{\de, s}(\f{M}{L})^{2(1+2\de(s-1))-s}\|u_0\|_{L^p}^2
\|u_{>M/2}(t, \cdot)\|_{L^2}^2,
\end{equation}
where the constantr $C_{\de, s}$ will depend on both $\de, s$ via the Sobolev
embedding estimate \eqref{eq:38}. Furthermore, by the log-convexity of $p\to
\|f\|_{L^p}$, we have $\|u_0\|_{L^p}\leq
\|u_0\|_{L^2}^{2/p}\|u_0\|_{L^\infty}^{1-2/p}$. In particular 
$\|u_0\|_{L^p}\leq \|u_0\|_{L^2\cap L^\infty}$, and we insert this  
in \eqref{eq:90}. 

Next, take $M=2^j L$, as this is somewhat more flexible for the forthcomming
induction argument.  We will show  the bound 
\eqref{eq:31} first for $0\leq t\leq 1$ and then, we will extend the result to
$t>1$. 
\subsubsection{$0\leq t\leq 1$}
We will show by induction that there exists a 
constant $C_0$, depending on $\de$ and $s$, so that 
\begin{equation}
\label{eq:01}
\|u_{>2^j L}(t)\|_{L^2}^2\leq (C_0 \max(1, \|u_0\|_{L^2\cap L^\infty}^2))^{j+1}
2^{-t(1-\de)(s-1) j^2}, \qq t\in[0,1].
\end{equation}
The first thing to observe is that for all $0<j\leq 5$, 
we have by the monotonicity of $t\to \|u(t)\|_{L^2}$, 
 $\|u_{>2^j L}(t)\|_{L^2}^2\leq \|u(t)\|_{L^2}^2\leq \|u_0\|_{L^2}^2$, whence  
 \eqref{eq:01} holds, 
as long as we select $C_0>2^{25(s-1)(1-\de)}$. \\
Thus, assuming the validity of \eqref{eq:01}  for some $j-1$, $j\geq 6$, we have by
\eqref{eq:90}
\begin{eqnarray*}
& & \p_t  \|u_{>2^j L} (t, \cdot)\|_{L^2}^2+ 2^{j s} 
\|u_{>2^jL}(t, \cdot)\|_{L^2}^2\leq C_{\de,s} 
 2^{j(2+2\de(s-1)-s)}\|u_0\|_{L^2\cap L^\infty}^2
\|u_{>2^{j-1}}(t, \cdot)\|_{L^2}^2\leq \\
& & 
\leq C_{\de,s} 2^{j(2+2\de(s-1)-s))}\|u_0\|_{L^2\cap L^\infty}^2 
(C_0 \max(1, \|u_0\|_{L^2\cap L^\infty}^2))^{j}
2^{-t(1-\de)(s-1) (j-1)^2},
\end{eqnarray*}
for every $0\leq t\leq 1$. 
Apply the Gronwall's inequality to the last equation. Note that to do that,  we 
have to take into account 
$$
\int_0^t e^{z(2^{js}-(j-1)^2(1-\de)(s-1))} dz\leq 2^{-js+1} 
e^{t2^{js}},
$$
 since $2^{js}>2 (j-1)^2$ for $j\geq 6, s>1$.  Thus,  
\begin{eqnarray*}
& &
\|u_{>2^j L} (t, \cdot)\|_{L^2}^2\leq \|P_{>2^j L} u_0\|_{L^2}^2 e^{-t2^{js}} + \\
& & + 
2 C_{\de,s} C_0^{j} \max(1, \|u_0\|_{L^2\cap L^\infty}^2)^{j+1} 2^{-2j(1-\de)(s-1)} 
2^{-t(1-\de)(s-1) (j-1)^2}.
\end{eqnarray*} 
The exponents that arise can be estimated  in the following  
straightforward manner. We have 
$e^{-t2^{js}}\leq  2^{-t(s-1)(1-\de)j^2}$ for all $j\geq 6, 1<s\leq 2, 1>\de>0$. Also,
since $t\in [0,1]$, we have $-2j(1-\de)(s-1)-t(1-\de)(s-1) (j-1)^2< -t(s-1)(1-\de)
j^2$. Thus, selecting \\ $C_0: C_0=4 C_{\de,s}+2^{25(s-1)(1-\de)}+2$ finishes the
proof of \eqref{eq:01}.
\label{sec:3.1.1}
\subsubsection{$t>1$} 
\label{sec:3.1.2}
The results of the previous case are easy to extend now to
the case $t>1$. Namely, we will show that there exists a constant $C_1$, so that 
\begin{equation}
\label{eq:02}
\|u_{>2^j L}(t)\|_{L^2}^2\leq (C_1 \max(1, \|u_0\|_{L^2\cap L^\infty}^2))^{j+1}
2^{-(1-\de)(s-1) j^2}, \qq t>1
\end{equation}
Again, the case of $j=0, \ldots, 5$ is easy to be verified by the monotonicity of
the $L^2$ norm. Assuming $j\geq 6$ and \eqref{eq:02} for all $t>1$ and some $j-1$,
we apply \eqref{eq:90}, where we insert the estimate \eqref{eq:02} 
for the term $u_{>2^{j-1}L}$. We get 
\begin{eqnarray*}
& & \p_t  \|u_{>2^j L} (t, \cdot)\|_{L^2}^2+ 2^{j s} 
\|u_{>2^jL}(t, \cdot)\|_{L^2}^2\leq \\
& & \leq  C_{\de,s} C_1^j 
2^{j(2+2\de(s-1)-s))} \max(1, \|u_0\|_{L^2\cap L^\infty}^2)^{j}
2^{-(1-\de)(s-1) (j-1)^2}.
\end{eqnarray*}
Apply the Gronwall's inequality in the interval $(1,t)$. 
\begin{eqnarray*}
& &
\|u_{>2^j L} (t, \cdot)\|_{L^2}^2\leq \|u_{>2^j L} (1, \cdot)\|_{L^2}^2 
e^{-(t-1)2^{js}}+ \\
& & +C_{\de,s} C_1^j \max(1, \|u_0\|_{L^2\cap L^\infty}^2)^{j+1}
2^{-2j(1-\de)(s-1)-(1-\de)(s-1) (j-1)^2}.
\end{eqnarray*}
However, inserting 
 the bound \eqref{eq:01} for $\|u_{>2^j L} (1, \cdot)\|_{L^2}^2$ and realizing that
 again  \\ 
 $-2j(1-\de)(s-1)-(1-\de)(s-1) (j-1)^2\leq -(1-\de)(s-1) j^2$, we have 
 for all $t>1$, 
\begin{eqnarray*}
& &
\|u_{>2^j L} (t, \cdot)\|_{L^2}^2\leq (C_0 \max(1, \|u_0\|_{L^2\cap
L^\infty}^2))^{j+1} 2^{-(1-\de)(s-1)j^2}+ \\
& & +
C_{\de,s} C_1^j \max(1, \|u_0\|_{L^2\cap L^\infty}^2)^{j+1} 
2^{-(1-\de)(s-1) j^2}\leq (C_1 \max(1, \|u_0\|_{L^2\cap L^\infty}^2))^{j+1}
2^{-(1-\de)(s-1) j^2},
\end{eqnarray*}
as long as $C_1=2C_0$. This concludes the proof of \eqref{eq:31}.

\subsection{The case $s>1+d/2$.} 
The proof for $s>1+d/2$ goes almost identically  to the case $1<s\leq 2$. 
Note that  
the monotonicity of $t\to \|u(t)\|_{L^p}, p>2$ is unavailable\footnote{Or at least,
we are not aware of such result.} in this context, but 
we still have that $t\to \|u(t)\|_{L^2}$ is decreasing and 
therefore by Lemma \ref{le:bern}
$$
\| u_{\leq M/2}\|_{L^\infty}\| + 
\|u_{\leq M}\|_{L^\infty} \leq C (M/L)^{d/2} \|u_0\|_{L^2},
$$
whence 
$$
\p_t  \|u_{>M} (t, \cdot)\|_{L^2}^2+ 2(M/L)^{s} 
\|u_{>M}(t, \cdot)\|_{L^2}^2 \leq C (M/L)^{1+d/2} 
\|u_{>M}\|_{L^2} \|u_{>M/2}\|_{L^2}  \|u_0\|_{L^2},
$$
whence
$$
\p_t  \|u_{>M} (t, \cdot)\|_{L^2}^2+ (M/L)^{s} 
\|u_{>M}(t, \cdot)\|_{L^2}^2 \leq C (M/L)^{2+d-s} \|u_{>M/2}\|_{L^2}^2
\|u_0\|_{L^2}^2.
$$
This is similar to \eqref{eq:90}, except for the power of $(M/L)$ on the right-hand
side. One can now perform an identical argument to show \eqref{eq:310}. This is done
by systematically replacing the factor $(1-\de)(s-1)$  by $s-1-d/2$, which is 
assumed to be positive.

\section{Estimates of the high-frequency tails for the 1 D
Kuramoto-Sivashinsky equation}
In this section, we  prove theorem \ref{theo:3}. The approach that we take is very
similar to the one in Section \ref{sec:10}, except that now because of the
destabilzing term $u_{xx}$, we do not have such a good control of $\|u(t)\|_{L^2}$. 
For the rest of the section, we will be proving \eqref{eq:42}.

We start as in Section \ref{sec:10} by taking the projection $P_{>M}$ in 
 \eqref{eq:2}, with $M>>L$.  After multiplication by $u$, integrating in $x$ and integration by
 parts, we obtain  
 \begin{eqnarray*}
 & & \p_t \f{1}{2}\|u_{>M}(t, \cdot)\|_{L^2}^2+
 \|\p_x^2 u_{>M}(t, \cdot)\|_{L^2}^2-\|\p_x u_{>M}(t, \cdot)\|_{L^2}^2\leq 
 |\int u_{>M} u u_x dx|
 \end{eqnarray*}
Now by the elementary properties of $P_{>M}$ in Section \ref{sec:four}, we have \\
$\|\p_x^2 u_{>M}(t, \cdot)\|_{L^2}^2\gtrsim (M/L)^{2} \|\p_x u_{>M}(t, \cdot)\|_{L^2}^2$ and
thus $\|\p_x^2 u_{>M}(t, \cdot)\|_{L^2}^2>>\|\p_x u_{>M}(t, \cdot)\|_{L^2}^2$. 
Moreover 
$ \|\p_x^2 u_{>M}(t)\|_{L^2}^2\gtrsim (M/L)^4 \|u_{>M}(t)\|_{L^2}^2$. 
On the other hand, following exactly the line of argument in Section \ref{sec:10}
\begin{eqnarray*}
 & &  |\int u_{>M} u u_x dx|\leq \f{1}{2} \|u_{>M}\|_{L^2}^2 \|\p_x u_{\leq
 M}\|_{L^\infty}+   \\ 
 & & + C (M/L) \|u_{>M}\|_{L^2}\|u_{>M/2}\|_{L^2} (\| u_{\leq M/2}\|_{L^\infty}\| + 
\|u_{\leq M}\|_{L^\infty} )
 \end{eqnarray*}
For  the second term on the right hand side, we further estimate via Cauchy-Schwartz 
\begin{eqnarray*}
 & &
(M/L) \|u_{>M}\|_{L^2}\|u_{>M/2}\|_{L^2} (\| u_{\leq M/2}\|_{L^\infty}\| + 
\|u_{\leq M}\|_{L^\infty} ) \leq \f{1}{4}(M/L)^4 \|u_{>M}\|_{L^2}^2+ \\
& &+ C (M/L)^{-2}  |u_{>M/2}\|_{L^2}^2(\| u_{\leq M/2}\|_{L^\infty}\| + 
\|u_{\leq M}\|_{L^\infty} )^2.
 \end{eqnarray*}
Putting all of these estimates together yields 
\begin{equation}
\label{eq:50}
\begin{array}{c}
\p_t  \|u_{>M}\|_{L^2}^2+
 2 (M/L)^4 \| u_{>M}\|_{L^2}^2\leq C \|u_{>M}\|_{L^2}^2 \|\p_x u_{\leq
 M}\|_{L^\infty}+ \\ 
+ C (M/L)^{-2} \|u_{>M/2}\|_{L^2}^2  (\| u_{\leq M/2}\|_{L^\infty} + 
\|u_{\leq M}\|_{L^\infty} )^2.
\end{array}
 \end{equation}
By Lemma \ref{le:bern}, 
$$
\|u_{\leq M}\|_{L^\infty}+\|u_{\leq M/2}\|_{L^\infty}\leq C (M/L)^{1/2} 
\supl_t \|u(t, \cdot)\|_{L^2}
$$
All in all, \eqref{eq:50}, together with the previous two observations implies 
\begin{eqnarray*}
& & \p_t  \|u_{>M}\|_{L^2}^2+
 2 (M/L)^4 \| u_{>M}\|_{L^2}^2 \leq  \\
 & & \leq C(M/L)^{3/2} 
\supl_s \|u(s, \cdot)\|_{L^2}\|u_{>M}\|_{L^2}^2 
+ C (M/L)^{-1} \sup_s \|u(s)\|_{L^2}^2 
 \|u_{>M/2}\|_{L^2}^2.
\end{eqnarray*}
Let $M=2^j L$ and denote $H=\sup_{s}\|u(s, \cdot)\|_{L^2}$. Fix an integer 
$j_0$, so that \\ 
$2^{5j_0}>100 \max(1, C^2) H^2$, where $C$ is the absolute constant appearing in
the last estimate.  In other words, our choice of $j_0$ is dictated by our 
need to ensure $2^{5j_0}>> H^2$. \\
 Denote 
 $I_j(t):=\|u_{> 2^j L}(t, \cdot)\|_{L^2}^2$. We have 
 \begin{equation}
 \label{eq:04}
 I_j'(t)+
 2^{4j+1} I_j(t) \leq 
 C 2^{3j/2} H I_j(t) + C 2^{-j} H^2 I_{j-1}(t). 
 \end{equation}
 Furthermore, since we are only interested in an estimate for $j\geq j_0$, it is
 easy to see that since $2^{-5j_0/2}H<<1$, 
 $$
 C 2^{3j/2} H I_j\leq C 2^{4j} 2^{-5j_0/2} H I_j< 2^{4j} I_j,
 $$
 which means that the first term on the right-hand side of \eqref{eq:04} may be
 absorbed on the left-hand side. Thus, for all $j\geq j_0$,  
 \begin{equation}
 \label{eq:05}
 I_j'(t)+
 2^{4j} I_j(t) \leq 
 C 2^{-j} H^2 I_{j-1}. 
 \end{equation}
We will apply the same idea as in the proof of \eqref{eq:31}. Namely, we  
run an induction argument based on \eqref{eq:05} for  $j\geq j_0$  
for a short period of time $0<t\leq 5/2$ and then we will extend to $t>5/2$. 
\subsection{$0\leq t\leq 5/2$}
We show that there exists an absolute constant $C_0$, so that for all 
$0<t<5/2$, and all $j\geq j_0$, 
\begin{equation}
\label{eq:06}
I_j(t)\leq C_0^{j+1} 2^{-t(j-j_0)^2} H^2.
\end{equation}
For $j=j_0$, the statement is obvious. Assuming the statement for some $j-1$, we
have by \eqref{eq:05}
\begin{eqnarray*}
& & I_j'(t)+
 2^{4j} I_j(t) \leq  C 2^{-j} H^2 C_0^j 2^{-t(j-1-j_0)^2} H^2
\end{eqnarray*}
Applying the Gronwall's inequality\footnote{Here again, we make use of the fact 
$\int_0^t exp(z(2^{4j}-(j-1-j_0)^2)) dz\leq 2^{-4j+1} exp(2^{4j}t)$, since 
$2^{4j}>2 (j-1-j_0)^2$, whenever $j\geq j_0+1$.} to the last inequality yields 
$$
I_j(t)\leq I_j(0)e^{-t 2^{4j}}+ C H^2 (2^{-5j_0} H^2) C_0^j 2^{-5(j-j_0)-t(j-1-j_0)^2}.
$$
Now, since $2^{4j}\geq (j-j_0)^2$ for $j\geq j_0$ and $I_j(0)\leq
\|u_0\|_{L^2}^2\leq H^2$, we have that 
$$I_j(0)e^{-t 2^{4j}}\leq 2^{-t(j-j_0)^2}H^2. 
$$
Next, since $C  (2^{-5j_0} H^2)<1$ and 
$-5(j-j_0)-t(j-1-j_0)^2\leq -t(j-j_0)^2$ (by $0<t\leq 5/2$), we conclude that 
$$
I_j(t)\leq C_0^{j+1} 2^{-t(j-j_0)^2} H^2, \qq t\in [0,5/2]
$$
whenever  $C_0\geq 2$. This concludes the proof of \eqref{eq:06}.
\subsection{$t > 5/2$} In this case, as in the Section \ref{sec:3.1.2}, we set
our induction argument with the hypothesis 
\begin{equation}
\label{eq:07}
I_j(t)\leq C_1^{j+1} 2^{-\f{5}{2}(j-j_0)^2} H^2.
\end{equation}
That is, we will show \eqref{eq:07} for all $j\geq j_0$ and for all $t>5/2$. 
We proceed as in Section \ref{sec:3.1.2}, namely we insert the induction hypothesis
in \eqref{eq:05} and then we run 
a Gronwall's argument  for the resulting inequality 
in the interval $[5/2,t]$. 

To give the proof in more detail, we start off with the observation that 
\eqref{eq:07} trivially holds with $j=j_0$. Assuming \eqref{eq:07} for some $j-1$,
we have by \eqref{eq:05}, 
\begin{eqnarray*}
& &  I_j'(t)+  2^{4j} I_j(t) \leq  C 2^{-j} H^2 
C_1^j 2^{-\f{5}{2}(j-1-j_0)^2} H^2. 
\end{eqnarray*}
By Gronwall's inequality, applied to the interval $[5/2,t]$, we have 
\begin{eqnarray*}
& &  I_j(t)\leq I_j(5/2) e^{(5/2-t) 2^{4j}}+ C 2^{-5j} H^2 C_1^j 
2^{-\f{5}{2}(j-1-j_0)^2} H^2\leq  \\ 
& & \leq I_j(5/2)+C C_1^j H^2 (2^{-5j_0}H^2) 
2^{-\f{5}{2}(2 (j-j_0)+(j-1-j_0)^2)}\leq \\
& & \leq C_0^{j+1} 2^{-\f{5}{2}(j-j_0)^2} H^2 + C C_1^j H^2 
2^{-\f{5}{2}(j-j_0)^2},
\end{eqnarray*}
where in the last inequality, we have used \eqref{eq:06} to estimate 
$I_j(5/2)$ and $2^{-5j_0} H^2<1$. Clearly, the last expression is estimated by 
$$
C_1^{j+1} 2^{-\f{5}{2}(j-j_0)^2} H^2,
$$
as claimed,  once we take
$C_1=2(C_0+C)$, where $C$ is an absolute constant appearing above. 

\section{Estimates of the higher Sobolev norms for the KSE}
In this section, we show how to make use of the Gevrey regularity estimates for
the solutions of KSE, provided by Theorem \ref{theo:3}, to provide effective
estimates on higher Sobolev norms. 
\subsection{Proof of Corollary \ref{cor:1}}
We actually show \eqref{eq:47}, which as we have showed implies \eqref{eq:43}. 
By the equivalence of the norms in \eqref{eq:020}, 
$$
\|u(t, \cdot)\|_{\dot{H}^s}\leq C^s [
\|u_{<C_0 H^{2/5}L}\|_{L^2} H^{2s/5}+ 
 \left(\suml_{j=0}^\infty (2^j C_0 H^{2/5})^{2s}
\|u_{\sim 2^j C_0 H^{2/5} L}\|_{L^2}^2\right)^{1/2}], 
$$
where $C$ is an absolute constant.  
For the first term, we have $\|u_{<C_0 H^{2/5}L}\|_{L^2}\leq H$. For 
the second term, we estimate by \eqref{eq:42}, 
$$
\supl_{\de\leq t} 
\|u_{\sim C_0 2^j H^{2/5}L}(t, \cdot)\|_{L^2}\leq 
C_1^j 2^{-\de j^2/2}H,
$$
which we insert in the sum above. We get 
\begin{eqnarray*}
& & \suml_{j=0}^\infty (2^j C_0 H^{2/5} )^{2s}
\|u_{\sim 2^j C_0 H^{2/5}}\|_{L^2}^2\leq C^s H^{4s/5+2} \suml_{j=0}^\infty 
C_1^{2 j} 2^{2sj-\de j^2} \leq C_{\de, s} H^{4s/5+2}. 
\end{eqnarray*} 
Taking square roots yields \eqref{eq:47}. 
\label{sec:lki}
\subsection{Proof of Corollary \ref{cor:2}.}
The proof of corollary \ref{cor:2} requires us to revisit the proof of Theorem
\ref{theo:3}. Namely, starting again with \eqref{eq:50}, we estimate this time 
(by Lemma \ref{le:bern})
$$
\|u_{\leq M}\|_{L^\infty}+\|u_{\leq M/2}\|_{L^\infty}\leq C_p (M/L)^{1/p} 
\supl_t \|u(t, \cdot)\|_{L^p}
$$
Thus, we get 
\begin{eqnarray*}
& &  \p_t  \|u_{>M}\|_{L^2}^2+
 2 (M/L)^4 \| u_{>M}\|_{L^2}^2 \leq  \\
 & & \leq C(M/L)^{1+1/p} 
\supl_s \|u(s, \cdot)\|_{L^p}\|u_{>M}\|_{L^2}^2 
+ C (M/L)^{-2+2/p} \sup_s \|u(s)\|_{L^p}^2 
 \|u_{>M/2}\|_{L^2}^2.
\end{eqnarray*}
Setting $M=2^j L$ and rewriting with $I_j(t)=\|u_{>2^j L}(t,
\cdot)\|_{L^2}^2$, we obtain the inequality 
\begin{equation}
\label{eq:021}
I_j'+ 2^{4j+1} I_j\leq C 2^{j(1+1/p)} K_p I_j + C 2^{j(-2+2/p)} K_p^2 I_{j-1}.
\end{equation}
Setting again $j_0: 2^{j_0(3-1/p)}=100\max(1, C^2) K_p$, we obtain that 
$$
C 2^{j(1+1/p)} K_p I_j\leq 2^{4j} I_j,
$$
and therefore one can absorb the first term on the right-hand side of
\eqref{eq:021}, as long as $j\geq j_0$. The result is 
$$
I_j'+ 2^{4j} I_j\leq C 2^{j(-2+2/p)} K_p^2 I_{j-1}.
$$
An induction argument similar to the one needed for the proof of  \eqref{eq:06} 
applies again. We get
\begin{equation}
\label{eq:022}
I_j(t)\leq C_0^{j+1} 2^{-t(j-j_0)^2}  H^2.
\end{equation}
for all $j\geq j_0$ and all\footnote{Note that in the previous argument, we have been
using $p=2$.} $t: 0<t<3-1/p$. 

In the case $t>3-1/p$, we apply an induction, similar to the one needed for the proof of
\eqref{eq:07}. We get for all $j\geq j_0$ and all $t>3-1/p$, 
$$
I_j(t)\leq C_1^{j+1}2^{-(3-1/p)(j-j_0)^2} H^2.
$$
Combining the two estimates yields the Gevrey bound 
\begin{equation}
\label{eq:023}
I_j(t)\leq C^{j+1} 2^{-\min(t, 3-1/p)(j-j_0)^2} H^2.
\end{equation}
Similarly to the proof of Corollary \ref{cor:1} (see Section \ref{sec:lki}), the Gevrey
estimate \eqref{eq:023} can be turned into estimates for higher Sobolev norms. Indeed, 
by \eqref{eq:023} and since $2^{j_0}\sim K_p^{1/(3-1/p)}$, we obtain  
$$
\supl_{\de\leq t} 
\|u_{\sim C_0 2^j K_p^{1/(3-1/p)}L}(t, \cdot)\|_{L^2}\leq 
C_1^j 2^{-\de j^2/2}H,
$$
whence \eqref{eq:012}.

\section{Characterization of the (eventual) blow-up time for the 2 D problem.}

Based on  classical results, we are assured that a solution is
classical up to the (eventual) blow-up time $T^*$. Thus, 
to show the characterization of
$T^*$ claimed in Theorem \ref{theo:2}, 
we proceed via a Gronwall inequality type
argument. 

Or first observation is that an integration in the $x$ variable in \eqref{eq:1}
yields  
\begin{equation}
\label{eq:5}
\p_t \int_{[-L,L]^2} -\phi(t,x) dx = 
\f{1}{2}\int_{[-L,L]^2}|\nabla\phi(t,x)|^2 dx.
\end{equation}
Next, we multiply \eqref{eq:1} by $\phi$ and integrate in the $x$ variable.
Keeping in mind that $\phi$ is real-valued and using integration by parts and
Cauchy-Schwartz, 
 we obtain 
\begin{eqnarray*}
& & \p_t \|\phi\|_{L^2}^2/2+ \|\De \phi\|_{L^2}^2-\|\nabla \phi\|_{L^2}^2 
\leq \f{1}{2} |\int \phi^2 (\De \phi)  dx|\leq \f{1}{2}\|\De
\phi\|_{L^2}\|\phi\|_{L^4}^2.
\end{eqnarray*}
At this point, we use the Sobolev embedding and Gagliardo-Nirenberg to estimate
$$
 \|\phi\|_{L^4}^2\leq C \|\phi\|_{H^{1/2}(-L,L)}^2\leq 
 C \|\nabla\phi\|_{L^2} \|\phi\|_{L^2}.
$$
All in all, after  
\begin{eqnarray*}
& & \p_t \|\phi\|_{L^2}^2/2+ \|\De \phi\|_{L^2}^2\leq 
\|\nabla \phi\|_{L^2}^2+ \f{1}{2} \|\De \phi\|_{L^2}^2 + C
\|\nabla\phi\|_{L^2}^2 \|\phi\|_{L^2}^2. 
\end{eqnarray*}
By the last inequality and  \eqref{eq:5}, it follows 
\begin{eqnarray*}
& & \p_t (\|\phi\|_{L^2}^2+1)\leq C (\|\phi\|_{L^2}^2+1)\|\nabla\phi\|_{L^2}^2=
C (\|\phi\|_{L^2}^2+1)\p_t \int(-\phi) dx.
\end{eqnarray*}
The Gronwall's inequality now implies 
$$
(\|\phi(t, \cdot)\|_{L^2}^2+1)\leq (\|\phi_0\|_{L^2}^2+1) 
exp(C \int(\phi_0(x)-\phi(t,x))dx)
$$
for some absolute constant $C$. 
The last inequality shows that the $\|\phi(t, \cdot)\|_{L^2}$ stays bounded
until either $\|\phi(t, \cdot)\|_{L^2}<\infty$ or  
$$
\int(-\phi(t,x))dx=\int(-\phi_0(x)dx+\f{1}{2}\int_0^T 
\|\nabla\phi(s, \cdot)\|_{L^2}^2 ds<\infty,
$$
which is satisfied, provided $\limsup_{T\to T^*} \int_0^T 
\|\nabla\phi(s, \cdot)\|_{L^2}^2 ds<\infty$. 

Analogously, one shows control over the higher order derivatives.  Let  
$\al=(\al_1, \al_2)$
to be a multindex in two variables, so that $\al_1, \al_2>2$. 
Then taking $\al$ derivatives of
\eqref{eq:1} and multiplying by $\p^\al \phi$ and integrating in $x$
yields\footnote{In what follows, for every integer $k$, 
 we use the notation $\al+k$ to denote the
multiindex $(\al_1+k, \al_2+k)$.} 
 \begin{equation}
\label{eq:11}
  \p_t \f{1}{2} \|\p^{\al}\phi\|_{L^2}^2+ \|\p^{\al+2} \phi\|_{L^2}^2- 
\|\nabla \p^{\al} \phi\|_{L^2}^2\leq \suml_{j=1}^2 |\int [\p^{\al} \phi] 
\p^{\al}[ \p_j 
\phi \p_j \phi ]dx|
\end{equation}
By integration by parts and Cauchy-Schwartz, we estimate the right-hand side 
$$
\suml_{j=1}^2 |\int [\p^{\al} \phi] \p^{\al} [ \p_j 
\phi \p_j \phi ]dx|\leq \f{1}{2} \|\p^{\al+2}\phi\|_{L^2}^2+ C 
\|\p^{\al-2}[\p\phi \p\phi]\|_{L^2}^2,
$$
where we have schematically denoted $\p \phi$ to stand for either derivative
$\p_1\phi, \p_2\phi$. 
We will need the following product estimate. 
\begin{lemma}
\label{le:1}
For every multindex $\al$ as above, there exists a constant $C_\al$, so that 
for every pair of   functions $u, v\in C^\infty_{per}([-L,L]^2)$, there is the estimate 
\begin{equation}
\label{eq:10}
\|\p^{\al-2}[\p u\p v]\|_{L^2}\leq C_\al(\||\nabla|^{|\al|} u\|_{L^2}\|\nabla
v\|_{L^2}+ \||\nabla|^{|\al|} v\|_{L^2}\|\nabla
u\|_{L^2})
\end{equation}
\end{lemma}
We postpone the proof of Lemma \ref{le:1}, so that we can finish our estimate showing control
of higher order Sobolev norms, with an {\it a priori} control of $\|\phi\|_{L^2}$. 
We have by  \eqref{eq:11}, and with the estimate of Lemma \ref{le:1}, 
we have established 
\begin{eqnarray*}
&  &
\p_t \f{1}{2} \|\p^{\al}\phi\|_{L^2}^2+ \|\p^{\al+2} \phi\|_{L^2}^2 \leq 
\f{1}{2} \|\p^{\al+2}\phi\|_{L^2}^2+C \|\p^\al \phi\|_{L^2}^2 \|\nabla\phi\|_{L^2}^2+
\|\p^\al\nabla\phi\|_{L^2}^2.
\end{eqnarray*}
For the last term on the right-hand side, we apply the Gagliardo-Nirenberg inequality 
$$
\|\p^\al\nabla\phi\|_{L^2}^2\leq \|\p^{\al+2}\phi\|_{L^2}^{2|\al|/(|\al|+1)}
\|\nabla\phi\|_{L^2}^{2/(|\al|+1)}\leq \f{1}{2} \|\p^{\al+2}\phi\|_{L^2}^2 
+ C_\al \|\nabla\phi\|_{L^2}^{2},
$$
where in the last inequality, we have used the Young's inequality $ab \leq a^p/p+b^q/q$,
for all $1<p,q<\infty: 1/p+1/q=1$. Putting these estimates together with \eqref{eq:5} yields  
$$
\p_t(\|\p^{\al}\phi\|_{L^2}^2+1)\leq
C_\al(\|\p^{\al}\phi\|_{L^2}^2+1)\|\nabla\phi\|_{L^2}^2=C_\al(\|\p^{\al}\phi\|_{L^2}^2+1) 
\p_t \int(-\phi) dx
$$
By Gronwall's, 
$$
\|\p^{\al}\phi(t)\|_{L^2}^2+1\leq (\|\p^{\al}\phi_0\|_{L^2}^2+1)exp(C_\al
\intl(\phi_0(x)-\phi(t,x)) dx),
$$
thus achieving the same control as before. 
\subsection{Proof of Lemma \ref{le:1}} 
By the alternative definition of $\|\cdot\|_{\dot{H}^s}$, it is enough
to show 
\begin{equation}
\label{eq:20}
(\suml_{j} 2^{2 j(|\al|-2)} \|P_{\sim 2^j L} (\p u \p v)\|_{L^2}^2)^{1/2}\leq C_\al
(\||\nabla|^{|\al|}  u\|_{L^2}\|\nabla
v\|_{L^2}+ \||\nabla|^{|\al|} v\|_{L^2}\|\nabla
u\|_{L^2})
\end{equation}
Furthermore, we have 
$$
P_{\sim 2^j L} (\p u \p v)=P_{\sim 2^j L} [\suml_{l_1} (\p P_{\sim 2^{l_1} L} 
u )(  \suml_{l_2} \p P_{\sim 2^{l_2} L} v)]
$$
Clearly, there are several cases to be considered, dependening on the relative strength
of $l_1$ to $j$. 
\subsubsection{$l_1>j+4$}
 Note that by support consideration,
$l_2>j+2$ and in fact $|l_2-l_1|\leq 2$. This is   because the product of two trig
polynomials, one of them of high degrees can be a low degree polynomial, if and only if
the two entries are of comparable degrees. We also observe that by the inclusion
$l^1\hookrightarrow l^2$, it is enough to replace the $l^2$ sum inthe left-hand side of
\eqref{eq:20}  with $l^1$ sum. 
Thus, the contribution of this piece (the so called ``high-high
interaction'')  is no more than 
\begin{eqnarray*}
& & \suml_{j} 2^{j(|\al|-2)} \|P_{\sim 2^j L} 
(\suml_{l_1>j+2} (\p P_{\sim 2^{l_1} L} u   \p v)\|_{L^2}= \\ 
& & \suml_{j} 2^{j(|\al|-2)} \|P_{\sim 2^j L} 
(\suml_{l_1>j+2} (\p P_{\sim 2^{l_1} L} u 
 \suml_{l_2: |l_2-l_1|\leq 2}  \p P_{\sim
2^{l_2} L}v)\|_{L^2}\leq  \\
& & C_\al \suml_{l_1} 2^{(|\al|-2)l_1} \|P_{\sim 2^{l_1}L} \p u\|_{L^\infty} 
\|P_{2^{l_1-2}L\leq \cdot\leq 2^{l_1+2}L} \p v\|_{L^2} \leq \\
& & C_\al (\suml_{l_1}  2^{2(|\al|-2)l_1}  
\|P_{\sim 2^{l_1}L} \p u\|_{L^\infty} ^2)^{1/2}  
(\suml_{l_1} \|P_{2^{l_1-2}L\leq \cdot\leq 2^{l_1+2}L} \p v\|_{L^2}^2)^{1/2}
\end{eqnarray*}
Clearly 
$$
(\suml_{l_1} \|P_{2^{l_1-2}L\leq \cdot\leq 2^{l_1+2}L} \p v\|_{L^2}^2)^{1/2}\lesssim
\|\p v\|_{L^2},
$$
while an application of Lemma \ref{le:bern} yields 
$$
\suml_{l_1}  2^{2(|\al|-2)l_1} 
\|P_{\sim 2^{l_1}L} \p u\|_{L^\infty} ^2\leq C \suml_{l_1}  2^{2(|\al|-2)l_1} 2^{2l_1}
\|P_{\sim 2^{l_1}L} \p u\|_{L^2} ^2\sim \||\nabla|^{|\al|} u\|_{L^2}^2.
$$
\subsubsection{$j-4<l_1<j+4$}
In that case one clearly has to have $l_2<j+6$. This
is simply because otherwise one will have a product of two polynomials - one of 
high degree ($l_2\geq j+6$) and one of low degree ($l_1<j+4$), and 
the resulting trig polynomial with degree $\sim 2^j L$, a contradiction. Thus, taking
into account Lemma \ref{le:bern}, the
contribution of this portion is less than 
\begin{eqnarray*}
& & (\suml_{j} 2^{2 j(|\al|-2)} \|P_{\sim 2^j L} (\p P_{2^{j-4}L <\cdot<2^{j+4}L} u)
(\p P_{<2^{j+6}L} v)\|_{L^2}^2)^{1/2}\lesssim   \\
& & (\suml_{j} 2^{2 j(|\al|-2)} \|\p P_{2^{j-4}L <\cdot<2^{j+4}L}
u\|_{L^\infty}^2)^{1/2}
\supl_m \|\p P_{<2^{m}L} v\|_{L^2}\lesssim \\
& &  (\suml_{j} 2^{2 j(|\al|-2)} 2^{4j} \|P_{2^{j-4}L <\cdot<2^{j+4}L}
u\|_{L^2}^2)^{1/2}
 \|\nabla v\|_{L^2}\leq  \||\nabla|^{|\al|} u\|_{L^2} 
\|\nabla v\|_{L^2}.
\end{eqnarray*}
\subsubsection{$l_1<j-4$} In this remaining case, similar Fourier support
considerations dictate that $l_2: |l_2-j|\leq 2$, and the estimate goes similarly to
the case $j-4<l_1<j+4$, just considered above.

\end{document}